\documentclass[11pt]{amsart}

\usepackage[hmargin=0.8in,twosideshift=0in,height=8.6in]{geometry}
\usepackage{amssymb,amsthm}
\usepackage{delarray,verbatim}
\usepackage{natbib}
\usepackage{ifpdf}
\ifpdf
\usepackage[pdftex]{graphicx}
\else
\usepackage[dvips]{graphicx}
\fi

\usepackage{ifpdf}
\usepackage{color}
\definecolor{webgreen}{rgb}{0,.5,0}
\definecolor{webbrown}{rgb}{.8,0,0}
\definecolor{emphcolor}{rgb}{0.95,0.95,0.95}

\usepackage{hyperref}
\hypersetup{%
          colorlinks=true,
          linkcolor=webbrown,
          filecolor=webbrown,
          citecolor=webgreen,
          breaklinks=true}
\ifpdf
\hypersetup{pdftex,
            pdfstartview=FitH, 
            bookmarksopen=true,
            bookmarksnumbered=true
}
\else
\hypersetup{dvips}
\fi

\newcommand {\ME}{\mathbb{E}^{x}}

\renewcommand{\S}{\mathcal{S}}

\numberwithin{equation}{section}
\newtheorem{proposition}{Proposition}[section]

\newtheorem{remark}{Remark}[section]
\newtheorem{lemma}[proposition]{Lemma}

\newcommand {\R}{\mathbb{R}}

\newcommand {\F}{\mathcal{F}}
\newcommand {\A}{\mathcal{A}}

\newcommand {\p}{\mathbb{P}}

\newcommand {\E}{\mathbb{E}}

\newcommand{\Norm}{N\left(\frac{x+\mu\Delta}{\sigma\sqrt{\Delta}}\right)}
\newcommand{\NormNeg}{N\left(\frac{-x+\mu\Delta}{\sigma\sqrt{\Delta}}\right)}
\newcommand{\Density}{\phi\left(\frac{x+\mu\Delta}{\sigma\sqrt{\Delta}}\right)}
\newcommand{\DensityNeg}{\phi\left(\frac{-x+\mu\Delta}{\sigma\sqrt{\Delta}}\right)}
\newcommand{\NormSqrt}{N\left(\frac{\sqrt{x}e^{-\rho\Delta}}{\sqrt{Q(\Delta)}}\right)}
\newcommand{\DensitySqrt}{\phi\left(\frac{\sqrt{x}e^{-\rho\Delta}}{\sqrt{Q(\Delta)}}\right)}

\title[]{A Unified Treatment of Dividend Payment Problems under Fixed Cost and Implementation Delays}
\author[]{Erhan Bayraktar }
\address[E. Bayraktar]{Department of
  Mathematics, University of Michigan, Ann Arbor, MI 48109}
\email{erhan@umich.edu}
\thanks{E. Bayraktar is supported in part by the National Science Foundation.}
\author[]{Masahiko Egami}
\address[M. Egami]{ Graduate School of Economics,
Kyoto University, Sakyo-Ku, Kyoto, 606-8501, Japan }
\email{egami@econ.kyoto-u.ac.jp}


\begin{document}

\begin{abstract}
In this paper we solve the dividend optimization problem for a
corporation or a financial institution when the managers of the
corporation are facing (regulatory) implementation delays. We
consider several cash reservoir models for the firm including  two
mean-reverting processes, Ornstein-Uhlenbeck and square-root
processes. Since the cashflow structure of different companies have different
qualitative behaviors it makes sense to use different diffusions to model them. We provide a
unified mathematical framework to analyze all these models and find the optimal barrier strategies. Our solution depends on a new characterization of
the value function for one-dimensional diffusions and provide
easily implementable algorithms to find the optimal control and
the value function.
\end{abstract}

\maketitle
\section{Introduction}

In this paper, we solve the dividend optimization problem for a
corporation or a financial institution. The corporation controls
the timing and the amount of dividends and the objective of the
corporation is to maximize the total discounted dividends paid out
to shareholders until the time of bankruptcy given that the
dividend payments are subject to regulatory delay. The payment of
a dividend is not automatic and payments can be made only after a
certain amount of time elapses.  The amount and the timing of
payment is decided by the company managers but these are subject
to the approval of the company's owners (shareholders) and maybe
also of debt holders and therefore it takes some time before the
dividends are paid. Recently, there have been other papers on
optimally controlling a state variable subject to implementation
delays in different modeling contexts, see e.g.
\cite{alvarez-keppo}, \cite{bar-ilan}, \cite{BE2006},
\cite{keppo-peura}, \cite{oksendal-delay-impulse} and
\cite{sub-jarrow}. Our methodolgy of solving this problem is in
the spirit of \cite{BE2006} and differs from the other papers
cited above as will be made clear below.

We model the problem of the corporation as an impulse control
problem and assume that when dividend is paid out, the firm has to
pay a fixed cost representing the resources it has to devote to
the distribution of dividends. This amount is independent of the
size of the dividend payment. Other papers modeling the dividend
payment problem as an impulse control problem are
\cite{Cadenillas}, \cite{shiryaev-jeanblanc} and \cite{paulsen07}. There are several
other papers which model the dividend payment problem as a
singular stochastic control problem by assuming that there is no
fixed cost at the time of dividend payment; see e.g. \cite{shiu-gerber-03},
\cite{shiu-gerber}, \cite{shiu-gerber-yang}, \cite{shiryaev-jeanblanc} and \cite{taksar}.

Applying an appropriate transformation to the value of a
particular control, we transform the problem into a non-linear
programming problem.  
 Using the
new characterization of the value function we give an easy to
implement algorithm to determine the optimal control and the value
function.
A secondary result of our paper are the
sufficient conditions under which the smooth fit holds (see Remark
4.1 and Proposition 4.1). 
In contrast, in the literature impulse control problems are solved first finding a classical solution to a system of quasi variational
inequalities. The optimal thresholds are determined using the so-called ``smooth fit principle" (by hypothesizing that the smooth fit holds).
(Once a classical solution to this system is determined
it can be shown to be equal to the value function by the so called \emph{verification lemma}.)
See e.g.
Bensoussan and Lions \cite{ben-lions} and {\O}ksendal and Sulem
\cite{oksendal-book-2}.


In this paper, the time horizon is the
time of ruin, and this makes the analysis more difficult from that of
\cite{BE2006}, which only considers infinite horizon problems. Since cashflow of different companies have different qualitative behavior, a manager needs a portfolio of tractable models to choose from. Here we consider four models for the aggregate
income/cash reservoir of the firm: i)\, Brownian motion with
drift, ii)\, Ornstein-Uhlenbeck, iii) Square-root process, iv)
Geometric Brownian motion. Most of the papers related to
stochastic impulse control, in order to obtain analytical
solutions, assume that the uncontrolled process is a Brownian
motion with drift. 
In addition to using
Brownian motion to model the cash
reservoir, we also propose two mean reverting processes as
possible modeling alternatives which is suggested by the Cash Flow
Hypothesis in Jensen \cite{Jensen}; see \cite{shiu-gerber-yang} for further motivation. 
On the other hand, geometric Brownian motion is used to model the firm value in the structural models in credit risk modeling. Our solution for the geometric Brownian motion model can also be interpreted as the optimal dividend distribution to the stockholders of a given company since geometric Brownian motion is frequently used to model the value of a company (for e.g. in the structural credit risk models) \cite{merton}; see \cite{shiu-gerber-03} for further motivation.
As far as we know, our paper is the first one that explicitly handles the dividend payment problem for
 for the square root process (with or without delays). 

The rest of the paper is organized as follows: In Section 2, we
present the models for the cash reservoir and state the dividend
payment problem. In Section 3, we provide a characterization of
the value function for a given threshold strategy. In Section 4,
we provide an easily implementable algorithm to find the optimal
threshold strategy. We also
provide theoretical justification for our algorithm in this
section (see e.g. Proposition~\ref{prop:main-proposition}). We then check that the models satisfy the sufficient assumptions of optimality in Section 4.3.
Finally, in Section 5 we present some numerical examples.

\section{Statement of The Problem}
Let $(\Omega, \F, \p)$ be a complete probability space with a
standard Brownian motion $W=\{W_t; t\geq 0\}$.  We model the
aggregate income process $X^0$ as either the Brownian motion
\begin{equation}\label{eq:process}
dX^0_t=\mu dt + \sigma dW_t, \quad X^0_0=x>0,
\end{equation}
for some constants $\mu, \sigma >0$; or the Ornstein Uhlenbeck
process
\begin{equation}\label{eq:OU}
dX^0_t=-\rho X^0_t dt + dW_t, \quad X^0_0=x>0,
\end{equation}
for some constant $\rho>0$, or the square root process
\begin{equation}\label{eq:CIR}
dX^0_t=(1-2 \rho X_t^0)dt+ 2 \sqrt{X^0_t} dW_t, \quad X^0_0=x>0.
\end{equation}
Note that if the initial condition of (\ref{eq:CIR}) is properly
chosen, then the solution of it is the square of the solution of
(\ref{eq:OU}). We will also consider the case when the aggregate
income process follows the geometric Brownian motion
\begin{equation}
dX_t^0=\mu X^{0}_t dt+ \sigma X^0_t dW_t, \quad X^{0}_0=x>0.
\end{equation}

The firm will pay dividends to its shareholders out of the
aggregate income process $X^0$ and the net holdings of the firm,
i.e. the net income process will be denoted by $X$. We assume that
the company pays out dividends to its shareholders in order to
maximize the expected value of discounted dividends paid out until
the time of ruin. There will be a fixed amount of transaction cost
for making a dividend payment. In this framework a dividend
payment scheme that a firm follows can be represented by a doubly
stochastic sequence
\begin{equation*}
\nu =(T_1, T_2,....T_i....; \xi_1,\xi_2,...\xi_i....),
\end{equation*}
where $0\leq T_1<T_2<....$ is an increasing sequence of
$\mathbb{F}$-stopping times such that $T_{i+1}-T_i \geq \Delta$,
and $\xi_1$, $\xi_2...$ are $\mathcal{F}_{(T_i+\Delta)-}$
measurable random variables representing the dividend amount paid out.
The firm decides to make dividend payments at (random) time $T_i$,
but it can not act until time $T_i+\Delta$ (where $\Delta \geq 0$
is a constant). It decides on the magnitude of the dividend amount
at $T_i+\Delta$ depending on the level of its revenues.
We will in particular consider benchmark strategies.  These strategies are determined by specifying two
numbers $0 \leq a<b$ as follows: At the time the aggregate
profit (or the firm value) hits a large enough level $b$, the
shareholders ask the firm to commit to making dividend payments
and reduce the level of net profits (or the firm value) to $a$. We denote
by $\mathcal{V}$ the set of strategies that fit into this
description. We will refer to them as the admissible strategies.

 The net income process follows
(until after the first dividend payment)
\begin{align} \label{eq:control}
\begin{cases}
dX_t = \mu(X_t)dt + \sigma(X_t)dW_t, \quad 0 \leq t<T_1+\Delta, \\
X_{T_{1}+\Delta}=X_{(T_1+\Delta)-}-\xi_1,
\end{cases}
\end{align}
for appropriate functions $\mu$ and $\sigma$ depending on which
case we are inspecting. For the first three cases we assume that
$0$ is the absorbing state and define $\tau_0$ (the time of ruin)
as :
\begin{equation*}
\tau_0\triangleq \inf \{t\geq 0: X_t\leq 0\}.
\end{equation*}
When the aggregate income process follows the geometric Brownian
motion, the time of ruin is defined as
\begin{equation}
\tau_{d} \triangleq \inf\{t \geq 0: X_t \leq d\},
\end{equation}
for some fixed $d>0$. The purpose of the firm is to maximize
expected value of the discounted dividend payments until the time
of ruin, i.e.,
\begin{align}\label{eq:J}
    J^\nu(x) \triangleq \mathbb{E}^x \left[\sum_{T_i+\Delta<\tau_0}e^{-\alpha
    (T_i+\Delta)} K(X_{(T_i+\Delta)-},X_{T_i+\Delta})\right],
\end{align}
over all the admissible strategies. We will assume that
\[
K(X_{(T_1+\Delta)-},X_{T_1+\Delta})=X_{(T_1+\Delta)-}-X_{T_1+\Delta}-\lambda,
\]
where
 $\lambda>0$ is a fee associated with a transaction. We could also
 consider
\[
K(X_{(T_1+\Delta)-},X_{T_1+\Delta})=k\cdot(X_{(T_1+\Delta)-}-X_{T_1+\Delta})-\lambda,
\]
for $k \in (0,1)$, in which $1-k$ can be considered as the tax
rate. This does not affect the analysis and therefore
we will focus on the case when $k=1$.

 Let us denote the value
function of this problem by
\begin{equation}\label{problem1}
v(x)\triangleq \sup_{\nu\in \mathcal{V}}J^{\nu}(x)=J^{\nu^*}(x).
\end{equation}

When $X^{0}$ is the Ornstein Uhlenbeck process, in addition to
considering the performance function in (\ref{eq:J}) we will also
consider the following performance
\begin{equation}\label{eq:new-penalty}
J^\nu(x) \triangleq \mathbb{E}^x
\left[\sum_{T_i+\Delta<\tau_0}e^{-\alpha
    (T_i+\Delta)} K(X_{(T_i+\Delta)-},X_{T_i+\Delta})-Pe^{-\alpha
    \tau_0}\right],
\end{equation}
for some constant $P>0$. The rationale for considering this
penalty function is to penalize declaring banktruptcy. As we shall
see if the purpose is to maximize the performance function in
(\ref{eq:J}), when $X^0$ follows an OU process, it is optimal to
declare bankruptcy when the aggregate income process reaches a
certain level. Therefore the OU process might be used to model the
income process of firms in distress. This might give an idea to
the creditors of how this type of a firm might behave. The extra
cost in (\ref{eq:new-penalty}) will, on the other hand, deter the
firms from declaring bankruptcy.

\section{Characterization of the Value
Function}\label{sec:characterization} In this section, we will
show that when we apply a suitable transformation (see \ref{defn:W}) to the value
function corresponding to a particular threshold strategy (that is
identified by a pair $(a,b)$), the transformed value function is
linear on $(F(0), F(b))$. This characterization will become important
in determining the optimal threshold strategy in the next section.
Equation (\ref{eq:J}) can be developed as
\begin{align*}
&J^{\nu}(x)=\ME\left[\sum_{T_i+\Delta<{\tau_0}}e^{-\alpha
    (T_i+\Delta)}K(X_{(T_i+\Delta)-},X_{T_i+\Delta})\right]\\
&= \ME\Big[1_{\{T_1+\Delta<\tau_0\}}\Big\{e^{-\alpha
(T_1+\Delta)}K(X_{(T_1+\Delta)-},X_{T_1+\Delta})\\
&\hspace{3.9cm}+e^{-\alpha(T_1+\Delta)}\sum_{i \geq 1:
T_i+\Delta<\tau_0,
}e^{-\alpha((T_{i+1}+\Delta)-(T_1+\Delta))}K(X_{(T_{i+1}+\Delta)-},X_{T_{i+1}+\Delta})\Big\}\Big]\\
&=\ME\Bigg[1_{\{T_1+\Delta<\tau_0\}}\Bigg\{e^{-\alpha
(T_1+\Delta)}K(X_{(T_1+\Delta)-},X_{T_1+\Delta})\\
&\hspace{4cm}+e^{-\alpha(T_1+\Delta)}\ME\left[\sum
e^{-\alpha((T_i+\Delta)\circ\theta(T_1+\Delta))}
K(X_{(T_{i+1}+\Delta)-},X_{T_{i+1}+\Delta})|\F_{T_1+\Delta}\right]\Bigg\}\Bigg]\\
&=\ME\left[1_{\{T_1+\Delta<\tau_0\}}e^{-\alpha
(T_1+\Delta)}\left\{K(X_{(T_1+\Delta)-},X_{T_1+\Delta})+\E^{X_{T_1+\Delta}}\sum_{T_i+\Delta<\tau_0}e^{-\alpha(T_i+\Delta)}K(X_{(T_i+\Delta)-},X_{T_i+\Delta})\right\}\right]
\\&=\ME\left[1_{\{T_1+\Delta<\tau_0\}}e^{-\alpha
(T_1+\Delta)}\left\{K(X_{(T_1+\Delta)-},X_{T_1+\Delta})+J^{\nu}(X_{T_1+\Delta})\right\}\right],
\end{align*}
where we used $T_{i+1}+\Delta=
(T_1+\Delta)+(T_i+\Delta)\circ\theta(T_1+\Delta)$ with the shift
operator $\theta(\cdot)$ in the second equality. 

\noindent Since $T_1=\tau_b$ with $\tau_b=\inf\{t\geq 0: X^0_t\geq
b\}$ and the post intervention point by
\begin{equation}
X_{T_1+\Delta}=X_{\tau_b+\Delta}=X_{(\tau_b+\Delta)-}-\xi_1\triangleq
a.
\end{equation}
It follows that
\begin{align}\label{eq:u2}
J^{\nu}(x)&=\ME\left[1_{\{\tau_b+\Delta<\tau_0\}}e^{-\alpha
(\tau_b+\Delta)}\left\{K(X_{(\tau_b+\Delta)-},
a)+J^{\nu}(a)\right\}\right]\nonumber\\
      &=\ME\left[\ME\left[1_{\{\tau_b+\Delta<\tau_0\}}e^{-\alpha
      (\tau_b+\Delta)}\left\{K(X_{(\tau_b+\Delta)-},
a)+J^{\nu}(a)\right\}\mid \F_{\tau_b}\right]\right]\nonumber\\
      &=\ME\left[1_{\{\tau_b<\tau_0\}}e^{-\alpha\tau_b}
      \E^{X_{\tau_b}}\left[1_{\{\Delta<\tau_0\}}e^{-\alpha\Delta}\left\{K(X_{\Delta-},
a)+J^{\nu}(a)\right\}\right] \right].\nonumber
\end{align}

Evaluating $J^{\nu}$ at $x=b$, we obtain
\begin{equation}\label{eq:x-at-b}
J^{\nu}(b)=\mathbb{E}^b\left[1_{\{\Delta<\tau_0\}}e^{-\alpha
\Delta}(K(X_{\Delta-},a)+J^{\nu}(a))\right],
\end{equation}
and evaluating $J^{\nu}$ at $x=0$ we get $J^{\nu}(0)=0$.

Now, we can write (\ref{eq:x-at-b}) as
\begin{equation}\label{u-cases}
J^{\nu}(x)=\begin{cases}\mathbb{E}^{x}\left[1_{\{\tau_b<\tau_0\}}e^{-\alpha
\tau_b} J^{\nu}(b)+ 1_{\{\tau_b>\tau_0\}}e^{-\alpha \tau_0}
J^{\nu}(0)\right] & 0\leq x \leq b \\
\mathbb{E}^x\left[1_{\{\Delta<\tau_0\}}e^{-\alpha \Delta
}(K(X_{\Delta},a)+J^{\nu}(a))\right] & x > b.
\end{cases}
\end{equation}
Similarly, if the performance function to minimize is the one
defined in (\ref{eq:new-penalty}), then we have
\begin{equation}\label{eq:jv}
J^{\nu}(x)=\begin{cases}\mathbb{E}^{x}\left[1_{\{\tau_b<\tau_0\}}e^{-\alpha
\tau_b} J^{\nu}(b)+ 1_{\{\tau_b>\tau_0\}}e^{-\alpha \tau_0}
J^{\nu}(0)\right], & 0\leq x \leq b, \\
\mathbb{E}^x\left[1_{\{\Delta<\tau_0\}}e^{-\alpha \Delta
}(K(X_{\Delta},a)+J^{\nu}(a))-P 1_{\{\Delta>\tau_0\}}e^{-\alpha
\tau_0 }\right] & x
> b.
\end{cases}
\end{equation}
In this case $J^{\nu}(0)=-P$ and
\[
J^{\nu}(b)=\mathbb{E}^b\left[1_{\{\Delta<\tau_0\}}e^{-\alpha
\Delta }(K(X_{\Delta},a)+J^{\nu}(a))-P
1_{\{\Delta>\tau_0\}}e^{-\alpha \tau_0 }\right].
\]
We will denote the infinitesimal generator of the process $X^0$ by
$\mathcal{A}$. Let us denote the increasing and decreasing
solutions of the second-order ordinary differential equation
$(\mathcal{A}-\alpha)u=0$ by $\psi(\cdot)$ and $\varphi(\cdot)$
respectively (these are uniquely determined up to a
multiplication). We can write
\begin{equation}\label{eq:laplace}
\ME[e^{-\alpha\tau_r}1_{\{\tau_r<\tau_l\}}]=\frac{\psi(l)\varphi(x)-\psi(x)\varphi(l)}
{\psi(l)\varphi(r)-\psi(r)\varphi(l)}, \, \,
\ME[e^{-\alpha\tau_l}1_{\{\tau_l<\tau_r\}}]=\frac{\psi(x)\varphi(r)-\psi(r)\varphi(x)}
{\psi(l)\varphi(r)-\psi(r)\varphi(l)},
\end{equation}
for $x\in[l,r]$ where $\tau_l\triangleq\inf\{t>0; X^0_t=l\}$ and
$\tau_r\triangleq\inf\{t>0; X^0_t=r\}$ (see e.g. Dayanik and
Karatzas \cite{DK2003}). Let us introduce the increasing function
\[
F(x) \triangleq \frac{\psi(x)}{\varphi(x)}.
\]
By defining
\begin{equation}\label{defn:W}
W \triangleq (J^{\nu}/\varphi)\circ F^{-1},
\end{equation}
on $x \in [0,b]$, using (\ref{eq:laplace}), equation
(\ref{u-cases}) on $ 0\leq x \leq b$ becomes
\begin{equation}\label{eq:linear}
W(F(x))=W(F(b))\frac{F(x)-F(0)}{F(b)-F(0)}+W(F(0))\frac{F(b)-F(x)}{F(b)-F(0)}
\quad 0\leq x \leq b
\end{equation}
which shows that the value function is \emph{linear} in the transformed space.
Next, we will compute
\begin{equation}\label{eq:defn-of-B}
B \triangleq \mathbb{E}^x\left[1_{\{\Delta<\tau_0\}}e^{-\alpha
\Delta }(K(X_{\Delta},a)+J^{\nu}(a))\right],
\end{equation}
in (\ref{u-cases}) for all the different models of aggregate
income process. (In the case of geometric Brownian motion we will replace $\tau_0$ by $\tau_d$ in \eqref{eq:defn-of-B}. Moreover the function $J^{\nu}(x)$ for this case is given by replacing $0$'s with $d$'s in \eqref{eq:jv}.) 
\subsection{Computation of $B$ in (\ref{eq:defn-of-B})}
\subsubsection{Ornstein-Uhlenbeck Process}
Let's first consider the case when $X^0$ is the Ornstein-Uhlenbeck
process given by (\ref{eq:OU}). Recall that $X^0_t$ can be written
as (can be derived using Theorem 4.6 of Karatzas and Shereve
\cite{kn:karat})
\begin{equation}
X^0_t=x e^{- \rho t}+ B_{Q(t)}, \quad \text{where} \quad
Q(t)=\frac{1-e^{-2 \rho t}}{2 \rho}, \quad \text{or}
\end{equation}
\begin{equation}\label{eq:mart}
e^{\rho t}X^0_t=x+ \tilde{B}_{\tilde{Q}(t)}, \quad \text{where}
\quad \tilde{Q}(t)=\frac{e^{2 \rho t}-1}{2 \rho},
\end{equation}
and $B$ and $\tilde{B}$ are Brownian motions. This implies that
the distribution of $X^0_t$ is $\text{Gsn} \left(xe^{-\alpha
t},Q(t)\right)$. (We use  $\text{Gsn}(a,b)$ to denote a Gaussian random variable with mean $a$ and variance $b$.) As a result of the representation in
(\ref{eq:mart})
\begin{equation}\label{eq:hit-time-zero}
\begin{split}
\mathbb{P}^x(\tau_0>\Delta)&=\mathbb{P}^{x}(\tau_0^{\tilde{B}}>\tilde{Q}(\Delta))=1-\frac{2}{\sqrt{2
\pi}} \int_{x/\sqrt{\tilde{Q}(\Delta)}}^{\infty}e^{-u^2/2}du
\\&=2
N\left(\frac{x}{\sqrt{\tilde{Q}(\Delta)}}\right)-1=2
N\left(\frac{x e^{-\rho \Delta}}{\sqrt{Q(\Delta)}}\right)-1,
\end{split}
\end{equation}
where $\tau_0^{\tilde{B}}$ is the first time the Brownian motion
$x+\tilde{B}$ hits zero. Here, we used the distribution of the
hitting times of Brownian motion (see page 96 of Karatzas and
Shreve). We also used the notation that
$N(x)=\int_{-\infty}^{x}1/(\sqrt{2 \pi})e^{-u^2/2}du$.

Let us try to identify the density function of $Y_{\Delta}
\triangleq  X^0_{\Delta}1_{\{\tau_{0}>\Delta\}}$. To this end we
first compute
\begin{equation}
\begin{split}
\mathbb{P}^{x}\{X^0_{\Delta} \geq y, \tau_0
>\Delta\}&=\mathbb{P}^{x}\{X^0_{\Delta} \geq y)-\mathbb{P}^{x}(X^0_{\Delta} \geq y, \tau_0 \leq
\Delta\}
\\&=\mathbb{P}^{x}\{X^0_{\Delta} \geq y\}-\mathbb{P}^x\{X^0_{\Delta} \leq -y, \tau_0 \leq
\Delta\}\\ &= \mathbb{P}^{x}\{X^0_{\Delta} \geq
y\}-\mathbb{P}^x\{X^0_{\Delta} \leq -y\}
\\ &=\frac{1}{\sqrt{2 \pi Q(\Delta)}} \left(\int_{y}^{\infty} \exp\left(-\frac{(u-x
e^{-\rho \Delta})^2}{2
Q(\Delta)}\right)du-\int_{-\infty}^{-y}\exp\left(-\frac{(u-x
e^{-\rho \Delta})^2}{2 Q(\Delta)}\right)du\right).
\end{split}
\end{equation}
Here, the second equality follows from the fact that OU process
satisfies a reflection principle around zero, and the third
inequality follows from the fact that $\{X^0_{\Delta} \leq -y\}
\supset\{\tau_0 \leq \Delta\}$ since $y>0$. The last line implies
that (after taking the derivative with respect to $y$ and flipping
the sign) the density of the random variable
$Y_{\Delta}=X^0_{\Delta}1_{\{\tau_{0}>\Delta\}}$ is given by
\begin{equation}\label{eq:density-of-stopped}
q(y)= \frac{1}{\sqrt{2 \pi Q(\Delta)}} \left(\exp\left(-\frac{(y-x
e^{-\rho \Delta})^2}{2 Q(\Delta)}\right)-\exp\left(-\frac{(y+x
e^{-\rho \Delta})^2}{2 Q(\Delta)}\right)\right), \quad y>0.
\end{equation}
Using (\ref{eq:hit-time-zero}) and (\ref{eq:density-of-stopped}),
we can write
\begin{equation}\label{eq:defn-B}
B= \mathbb{E}^x\left[1_{\{\Delta<\tau_0\}}e^{-\alpha \Delta
}(K(X_{\Delta},a)+J^{\nu}(a))\right]=e^{-\alpha \Delta
}\left(\left(2 N\left(\frac{x e^{-\rho
\Delta}}{\sqrt{Q(\Delta)}}\right)-1\right)(J^{\nu}(a)-a-\lambda)+A\right),
\end{equation}
in which
\[
A\triangleq\int_{0}^{\infty}y q(y)dy.
\]
Since
\[
\int_{0}^{x e^{-\rho \Delta}} y \exp\left(-\frac{(y-x e^{-\rho
\Delta})^2}{2 Q(\Delta)}\right)dy=-\int_{-x e^{-\rho \Delta}}^{0}y
 \exp\left(-\frac{(y+x e^{-\rho \Delta})^2}{2 Q(\Delta)}\right)dy,
\]
we can write $A$ as
\[
A=\frac{1}{\sqrt{2 \pi Q(\Delta)}}\left(\int_{x e^{-\rho
\Delta}}^{\infty} y \exp\left(-\frac{(y-x e^{-\rho \Delta})^2}{2
Q(\Delta)}\right)dy -\int_{-x e^{-\rho \Delta}}^{\infty} y
\exp\left(-\frac{(y+x e^{-\rho \Delta})^2}{2
Q(\Delta)}\right)dy\right).
\]
For any $\mu \in \mathbb{R}$ and $\sigma>0$ we have that
\[
\int_{\mu}^{\infty}\frac{x}{\sqrt{2 \pi
\sigma^2}}\exp\left(-\frac{(x-\mu)^2}{2
\sigma^2}\right)dx=\frac{\sigma}{\sqrt{2 \pi}}+\frac{\mu}{2}.
\]
As a result
\begin{equation}
A=x e^{-\rho \Delta}.
\end{equation}
Observe from the above calculations that
$X^0_{\Delta}1_{\{\tau_{0}>\Delta\}}$ and $X^0_{\Delta}$ have the
same expectation.

We will also compute the quantity
\begin{equation}
\tilde{B} \triangleq
\mathbb{E}^x\left[1_{\{\Delta<\tau_0\}}e^{-\alpha \Delta
}(K(X_{\Delta},a)+J^{\nu}(a))-P 1_{\{\Delta>\tau_0\}}e^{-\alpha
\tau_0}\right],
\end{equation}
for this case. Using the density of the hitting time of 0, which
can be derived by differentiating (\ref{eq:hit-time-zero}) we can
write
\begin{equation}\label{eq:tilde-B}
\tilde{B}=B-P \int_{0}^{\Delta}e^{-\alpha t} \frac{x}{\sqrt{2
\pi}}\left(\frac{\rho}{\sinh(\rho
t)}\right)^{\frac{3}{2}}\exp\left(- \frac{\rho x^2 e^{-\rho t }}{2
\sinh(\rho t)}+\frac{\rho t}{2}\right)dt.
\end{equation}
There is not explicit expression available for the integral term
(even in terms of special functions, except when $\Delta=\infty$,
see e.g. \cite{salminen} and \cite{darling}, in which case this
integral is the Laplace transform of the distribution of $\tau_0$)
but the NIntegrate function of Mathematica is able to evaluate it
with a very high numerical precision.

\begin{remark}
We can compute $B$ in (\ref{eq:defn-of-B}) explicitly even for the
cases when $X^0$ follows
\begin{equation}
dX^{0}_t=(\phi-\rho X^0_t)dt+\sigma dW_t, \quad X_0=x>0,
\end{equation}
for $\phi,\sigma>0$ by using the Strong Markov property to compute
\[
\E^x[X_{\Delta}1_{\{\Delta<\tau_0\}}]=\E^x[X_{\Delta}]-\E^x[1_{\{\Delta
\geq \tau_0\}}X_{\Delta}].
\]
The Strong Markov property is used to compute
\[
\begin{split}
\E^x[1_{\{\Delta \geq \tau_0\}}X_{\Delta}]&=\E^x[1_{\{\Delta \geq
\tau_0\}}\E^x[X_{\Delta}|\mathcal{F}_{\tau_0}]] =\E^x[1_{\{\Delta
\geq \tau_0\}}\E^0[X_{\Delta-\tau_0}]]
\\&=\int_{0}^{\Delta}f(u)E^{0}[X_{\Delta-u}]=\phi
\int_{0}^{\Delta}f(u)(1-\exp(-\rho(\Delta-u)))du,
\end{split}
\]
where $f$ is the density function of $\tau_0$. Several
representations for $f$ are available, see for e.g. \cite{alili}.

\end{remark}

\subsubsection{Square-root Process}
To evaluate $B$ in (\ref{eq:defn-of-B}) when the aggregate income
process is modeled by the square root process in (\ref{eq:CIR}) we
need to compute
\begin{equation}
C\triangleq \int_{0}^{\infty}y^{2}\tilde{q}(y)dy,
\end{equation}
in which $\tilde{q}(y)$ is equal to the $q(y)$ in
(\ref{eq:density-of-stopped}) if $x$ is replaced by $\sqrt{x}$.
This follows because if (\ref{eq:OU}) is started from $\sqrt{x}$,
then the solution of it is the square root of the solution of
(\ref{eq:CIR}). Let us first evaluate {\small
\begin{equation}\label{eq:eval-C}
\begin{split}
&\frac{1}{\sqrt{2 \pi
Q(\Delta)}}\int_0^{\infty}y^{2}\exp\left(-\frac{(y-\sqrt{x}e^{-\rho
\Delta})^2}{2 Q(\Delta)}\right)dy = \frac{1}{\sqrt{2
\pi}}\int_{-\frac{\sqrt{x}e^{-\rho
\Delta}}{\sqrt{Q(\Delta)}}}^{\infty}(y \sqrt{Q(\Delta)}+\sqrt{x}
e^{-\rho \Delta})^2\exp\left(-\frac{y^2}{2}\right)dy
\\ & =\frac{Q(\Delta)}{\sqrt{2 \pi}} \int_{-\frac{\sqrt{x}e^{-\rho
\Delta}}{\sqrt{Q(\Delta)}}}^{\infty}y^2
\exp\left(-\frac{y^2}{2}\right)dy +\frac{2
\sqrt{Q(\Delta)}}{\sqrt{2 \pi}} \sqrt{x} e^{-\rho \Delta
}\int_{-\frac{\sqrt{x}e^{-\rho
\Delta}}{\sqrt{Q(\Delta)}}}^{\infty}y
\exp\left(-\frac{y^2}{2}\right)dy+ x e^{-2 \rho \Delta}
N\left(\frac{\sqrt{x}e^{-\rho \Delta}}{\sqrt{Q(\Delta)}}\right)
\\ &= \frac{Q(\Delta)}{\sqrt{2 \pi}} \int_{-\frac{\sqrt{x}e^{-\rho
\Delta}}{\sqrt{Q(\Delta)}}}^{\infty}y^2
\exp\left(-\frac{y^2}{2}\right)dy +\frac{2
\sqrt{Q(\Delta)}}{\sqrt{2 \pi}} \sqrt{x} e^{-\rho \Delta }
\exp\left(-\frac{x e^{-2 \rho \Delta}}{2 Q(\Delta)}\right)+x e^{-2
\rho \Delta} N\left(\frac{\sqrt{x}e^{-\rho
\Delta}}{\sqrt{Q(\Delta)}}\right).
\end{split}
\end{equation}}
Since $\frac{1}{\sqrt{2 \pi
Q(\Delta)}}\int_0^{\infty}y^{2}\exp\left(-\frac{(y+\sqrt{x}e^{-\rho
\Delta})^2}{2 Q(\Delta)}\right)dy$ can be obtained by flipping the
sign in front of $\sqrt{x}$ in (\ref{eq:eval-C}), the computation
of $C$ will follow. We also have that
\begin{equation}\label{eq:int-second}
\begin{split}
&\frac{Q(\Delta)}{\sqrt{2 \pi}} \int_{-\frac{\sqrt{x}e^{-\rho
\Delta}}{\sqrt{Q(\Delta)}}}^{\infty}y^2
\exp\left(-\frac{y^2}{2}\right)dy = \frac{Q(\Delta)}{
2}+\frac{Q(\Delta)}{\sqrt{2 \pi}} \int_{-\frac{\sqrt{x}e^{-\rho
\Delta}}{\sqrt{Q(\Delta)}}}^{0}y^2
\exp\left(-\frac{y^2}{2}\right)dy
\\&= \frac{Q(\Delta)}{ 2}+ Q(\Delta)N\left(\frac{\sqrt{x}e^{-\rho
\Delta}}{\sqrt{Q(\Delta)}}\right)-\frac{\sqrt{x
Q(\Delta)}}{\sqrt{2 \pi}}e^{-\rho \Delta}\exp\left(-\frac{x e^{-2
\rho \Delta}}{2 Q(\Delta)}\right).
\end{split}
\end{equation}
Form (\ref{eq:eval-C}) and (\ref{eq:int-second}), we can evaluate
$C$ as
\begin{equation}
C= (x e^{-2 \rho \Delta}+ Q(\Delta))\left[2 N \left(\frac{\sqrt{x}
e^{-\rho \Delta}}{\sqrt{Q(\Delta)}}\right)-1\right]+ \frac{\sqrt{2
Q(\Delta)x}}{\sqrt{\pi}}e^{-\rho \Delta} \exp\left(-\frac{x e^{-2
\rho \Delta}}{Q(\Delta)}\right).
\end{equation}
When $X^0$ is the square root process then $B$ defined in
(\ref{eq:defn-of-B}) equals
\begin{equation}
B=e^{-\alpha \Delta }\left(\left(2 N\left(\frac{\sqrt{x} e^{-\rho
\Delta}}{\sqrt{Q(\Delta)}}\right)-1\right)(J^{\nu}(a)-a-\lambda)+C\right).
\end{equation}
\subsubsection{Brownian Motion with Drift}
Similarly, using reflection principle, Girsanov's Theorem and the
spatial homogeneity of Brownian motion we will obtain $B$ in
(\ref{eq:defn-B}) when $X^0$ is a Brownian motion given by
(\ref{eq:process}). We will first need the following lemma, which
is Corollary B.3.4 in \cite{musiela-rutkowski}.
\begin{lemma}\label{lemma:musiela}
Let $Y_t=\sigma W_t+\mu t$ and $m_t^Y=\min_{u \in [0,t]}Y_u$. Then
\begin{equation}
\mathbb{P}\{Y_t \geq y, m_t^{Y} \geq m\}=N\left(\frac{-y+\mu
t}{\sigma \sqrt{t}}\right)-e^{2 \mu m/\sigma^2}
N\left(\frac{2m-y-\mu t}{\sigma \sqrt{t}}\right),
\end{equation}
for every $m \leq 0$ and $y \geq m$.
\end{lemma}

We can write $B$ as
\begin{equation}
B=\E^{x+\theta}\left[1_{\{\Delta<\tau_{\theta}\}} X^0_{\Delta}
\right],
\end{equation}
in which $\theta \triangleq J^{\nu}(a)-a-\lambda$, which follows
from the spatial homogeneity of Brownian motion. Note that for any
$y>\theta$
\[
Z  \triangleq 1_{\{\Delta<\tau_{\theta}\}} X^0_{\Delta} \geq y
\quad \text{iff} \quad X_{\Delta} \geq y, \, m^{X^0}_{\Delta} \geq
\theta.
\]
We will find the probability density function of Z. Let us first
define
\[
Y_t \triangleq X_t-(\theta+x),
\]
which implies that $m_t^{Y}=m_t^{X^0}-(\theta+x)$. With this new
definition
\begin{equation}
\mathbb{P}\{Z \geq y\}=\mathbb{P}\{Y_{\Delta} \geq y-(x+\theta),
m_{\Delta}^Y \geq -x\}= N\left(\frac{-y+x+\theta+\mu
\Delta}{\sigma \sqrt{\Delta}}\right)-e^{-2 \mu x/\sigma^2}
N\left(\frac{-x-y+\theta+ \mu \Delta}{\sigma
\sqrt{\Delta}}\right).
\end{equation}
Here, the second equality follows from Lemma~\ref{lemma:musiela}.
Now, the density of the random variable $Z$ is easy to calculate
and using that we can compute $B$ by calculating the expectation
of $Z$ and get
\begin{equation}
\begin{split}
B&= e^{-\alpha \Delta}\Bigg((x+\mu
\Delta+J^{\nu}(a)-a-\lambda)N\left(\frac{x+\mu \Delta }{\sigma
\sqrt{\Delta}}\right) +\frac{\sigma \sqrt{\Delta}}{\sqrt{2\pi}}
 \exp\left(-\frac{1}{2}\frac{(x+\mu \Delta)^2 }{\sigma^2
\Delta}\right) \\&-e^{-2 \mu x/\sigma^2}\bigg((-x+\mu \Delta+
J^{\nu}(a)-a-\lambda) N\left(\frac{-x+\mu \Delta }{\sigma
\sqrt{\Delta}}\right)+  \frac{\sigma \sqrt{\Delta}}{\sqrt{2 \pi}}
\exp\left(-\frac{1}{2}\frac{(-x+\mu \Delta)^2 }{\sigma^2
\Delta}\right) \bigg) \Bigg).
\end{split}
\end{equation}

\subsubsection{Geometric Brownian Motion}

We will use the down and out European call option price, see e.g.
\cite{rubinstein} , (when we take the strike price to be zero) to
evaluate
\begin{equation}
\E^x\left[1_{\{\tau_d>\Delta\}} X_{\Delta} \right]=e^{\mu \Delta}
x \left(N(d_1)- \left(\frac{d}{x}\right)^{1+2 \mu/\sigma^2}
N(-d_2)\right),
\end{equation}
in which
\begin{equation}
\begin{split}
d_1=\frac{\log\frac{x}{d}+(\mu-\frac{1}{2}\sigma^2)\Delta}{\sigma
\sqrt{\Delta}}  \quad d_2=
\frac{\log\frac{x}{d}-(\mu-\frac{1}{2}\sigma^2)\Delta}{\sigma
\sqrt{\Delta}}.
\end{split}
\end{equation}
In order to calculate $B$ we also need to compute
$\p\{\tau_d>\Delta\}$. In fact
\[
\p^x\{\tau_d>\Delta\}=\p^x\{\tau^{\tilde{B}}_{\tilde{d}}>\Delta\},
\]
in which $\tau^{\tilde{B}}_{\tilde{d}}$ is the hitting time of
$\tilde{d}<0$ by the Brownian motion $\tilde{B}\triangleq \gamma
t+ \sigma B_t$, where
\[
\tilde{d} \triangleq \log \left(\frac{d}{x}\right), \quad
\gamma=\mu-\frac{1}{2}\sigma^2.
\]
Using the hitting time distribution for Brownian motion (with
drift), (which can be obtained from Lemma~\ref{lemma:musiela}), we
deduce
 \begin{equation}
 \p^x\{\tau_d>\Delta\}=\p^x\{m^{\tilde{B}}_{\Delta} \geq \tilde{d}\}=N\left(\frac{-\tilde{d}+\gamma t
 }{\sigma \sqrt{\Delta}}\right)-\exp\left(\frac{2 \gamma \tilde{d}}
 {\sigma^2}\right)N\left(\frac{\tilde{d}+\gamma t}{\sigma \sqrt{\Delta}}\right).
 \end{equation}
 Therefore, $B$ can be written as
 \begin{equation}
 B=e^{-\alpha \Delta} \left[e^{\mu \Delta} x \left(N(d_1)- \left(\frac{d}{x}\right)^{1+2 \mu/\sigma^2} N(-d_2)\right)
+(J^{\nu}(a)-a-\lambda) \p^x\{\tau_d>\Delta\}\right].
 \end{equation}

\section{An Efficient Algorithm to Calculate the Value Function}
\subsection{Increasing and Decreasing Solutions of $(\mathcal{A}-\alpha) u=0$}
When $X^0$ is the Brownian motion in (\ref{eq:process}), then the
increasing and decreasing solutions of $(\mathcal{A}-\alpha)u=0$
are
\begin{equation}
\psi(x)=e^{D_1x} \quad \text{and} \quad \varphi(x)=e^{D_2x}
\end{equation}
in which
\begin{equation*}
D_1=\frac{-\mu+\sqrt{\mu^2+2\alpha\sigma^2}}{\sigma^2}\quad\text{and}\quad
D_2=\frac{-\mu-\sqrt{\mu^2+2\alpha\sigma^2}}{\sigma^2}.
\end{equation*}
 When $X^0$ is the
Ornstein-Uhlenbeck process in (\ref{eq:OU}), then
\begin{equation}
\psi(x)=\exp\left(\frac{\rho x^2}{2}\right)D_{-\alpha/\rho}(-x
\sqrt{2 \rho}), \quad \varphi(x)=\exp\left(\frac{\rho
x^2}{2}\right)D_{-\alpha/\rho}(x \sqrt{2 \rho}), \quad x \in
\mathbb{R},
\end{equation}
where $D_{\nu}(\cdot)$ is the parabolic cylinder function given in
the Appendices 1.14 and 2.9 in \cite{salminen} which is defined as
\[
D_{\nu}(x)\triangleq
2^{-\nu/2}e^{-x^2/4}H_{v}\left(\frac{x}{\sqrt{2}}\right),\quad x
\in \mathbb{R},
\]
in which $H_{\nu}$ is the Hermite polynomial of order $\nu$, which
has the integral representation (see e.g. \cite{lebedev})
\[
H_{\nu}(x)=\frac{1}{\Gamma(-\nu)}\int_{0}^{\infty}\exp\left(-t^2-2tx\right)t^{-\nu-1}dt,
\quad \text{Re}(\nu)<0.
\]
On the other hand, when $X^0$ is the square root process whose
dynamics follows (\ref{eq:CIR}), then
\begin{equation}
\psi(x)=x^{-1/4}\exp\left(\frac{\rho
x}{2}\right)M_{-\frac{\alpha}{2
\rho}+\frac{1}{4},-\frac{1}{4}}(\rho x), \quad
\varphi(x)=x^{-1/4}\exp\left(\frac{\rho
x}{2}\right)W_{-\frac{\alpha}{2
\rho}+\frac{1}{4},-\frac{1}{4}}(\rho x),
\end{equation}
in which $W_{-\frac{\alpha}{2 \rho}+\frac{1}{4},-\frac{1}{4}}$ and
$M_{-\frac{\alpha}{2 \rho}+\frac{1}{4},-\frac{1}{4}}$ are
Whittaker functions (see e.g. Appendix 2.10 of \cite{salminen}).
These functions satisfy
\begin{equation}
\begin{split}
W_{-\frac{\alpha}{2
\rho}+\frac{1}{4},-\frac{1}{4}}\left(\frac{x^2}{2}\right)&=2^{\frac{\alpha}{2
\rho}-\frac{1}{4}}\sqrt{x}D_{-\alpha/\rho}(x), \quad x \geq 0,
\\M_{-\frac{\alpha}{2
\rho}+\frac{1}{4},-\frac{1}{4}}\left(\frac{x^2}{2}\right)&=\frac{\Gamma((1+\alpha/\rho)/2)}{2
\sqrt{\pi}}\sqrt{x}(D_{-\alpha/\rho}(-x)-D_{-\alpha/\rho}(x)),
\quad x \geq 0,
\end{split}
\end{equation}
in which $\Gamma$ stands for the Gamma function
$
\Gamma(x)=\int_{0}^{\infty}u^{x-1}e^{-u}du.
$
When, $X^0$ is the geometric Brownian motion, then
\begin{equation}
\psi(x)=x^{\sqrt{\kappa^2+\frac{2 \alpha}{\sigma^2}}-\kappa} \quad
\varphi(x)=x^{-\sqrt{\kappa^2+\frac{2 \alpha}{\sigma^2}}-\kappa},
\end{equation}
in which $\kappa=\mu/\sigma^2-1/2$.

\subsection{An Algorithm to Find the Optimal Control}
In this section we will describe a numerical algorithm to find the
value function. First we will introduce some notation that we
facilitate our description.
\begin{equation}\label{defn:r-h}
e^{\alpha \Delta} B =: r(x;a)+ h(x) J^{\nu}(a),
\end{equation}
where $B$ is as in (\ref{eq:defn-of-B}). We transform $r$ and $h$
by
\begin{equation}\label{eq:transform-r-h}
R(\cdot; a)\triangleq \frac{r(F^{-1}(\cdot),
a)}{\varphi(F^{-1}(\cdot))}, \quad H(\cdot) \triangleq
\frac{h(F^{-1}(\cdot))}{\varphi(F^{-1}(\cdot))}.
\end{equation}
Note that $r(a;a)<0$ and that  $\sup_{x}r(x;a)>0$ in all the cases considered above (see Section~\ref{sec:satisfaction}). 
\noindent \underline{First
stage}: For a given pair $(a,b) \in \mathbb{R}_{+}^{2}$ we will
determine $W$ in (\ref{defn:W}) using the linear characterization
in (\ref{eq:linear}). On $[F(0),F(b)]$ we will find the line
$W(y)=\beta y+\xi$ that passes through the point
$\left(F(0),-\frac{P}{\varphi(0)}\right)$, i.e.,
\begin{equation}
\xi=-\beta F(0)-\frac{P}{\varphi(0)},
\end{equation}
and satisfies
\begin{equation}\label{eq:continuous-fit}
\beta F(b)+ \xi= e^{-\alpha \Delta} \left[R(F(b);a)+H(F(b))
\varphi(a) (\beta F(a)+\xi) \right].
\end{equation}
$P=0$ when we consider (\ref{u-cases}). The slope $\beta$ can be
determined as
\begin{equation}\label{eq:slope}
\beta=\frac{\frac{P}{\varphi(0)}+e^{-\alpha
\Delta}\left[R(F(b);a)-H(F(b)) \varphi(a)
\frac{P}{\varphi(0)}\right]}{F(b)-F(a)-e^{-\alpha
\Delta}(F(a)-F(0))H(F(b))\varphi(a)}\, .
\end{equation}
Now the function $J^{\nu}(x)$ in \eqref{eq:jv} can be written as
\begin{equation}\label{eq:perform-a-b}
J^{\nu}(x)=
\begin{cases}
\beta \psi(x)+ \xi \varphi(x), & 0 \leq x \leq b,
\\ e^{-\alpha \Delta} (r(x;a)+ h(x) J^{\nu}(a)), & x\geq b.
\end{cases}
\end{equation}
Note that $(\mathcal{A}-\alpha) J^{\nu}(x)=0$ for $x<b$.

\noindent \underline{Second stage}: Let us fix $a$ and treat
$\beta$ as a function of $b$ parametrized by $a$. We will maximize
the function $\beta$ in (\ref{eq:slope}). Taking the derivative of
(\ref{eq:continuous-fit}) and evaluating at $\beta_b=0$ we obtain
\begin{equation}\label{eq:find-b-given-a}
\beta =e^{-\alpha \Delta} \left[\frac{\partial}{\partial y}R(y;a)
\bigg|_{y=F(b)} + H'(F(b)) \varphi(a) \left(\beta \cdot (F(a)-
F(0))-\frac{P}{\varphi(0)}\right)\right],
\end{equation}
in which $\beta$ is as in (\ref{eq:slope}). To find the optimal
$b$ given $a$ we solve the non-linear and implicit equation
(\ref{eq:find-b-given-a}).
\begin{remark}\label{remark:smooth-fit}
On $y \geq F(b)$, the function $W$ is given by
\begin{equation}\label{eq:continuous2-fit}
W(y)=e^{-\alpha \Delta}\left[W(F(a))\varphi(a)H(y)+R(y;a)\right].
\end{equation}
The right derivative of $W$ at $F(b)$ is
\begin{equation}\label{eq:smooth-fit}
W'(F(b))=e^{-\alpha
\Delta}\left(W(F(a))\varphi(a)H'(F(b))+\frac{\partial}{\partial y
}R(y;a)\bigg|_{y=F(b)}\right)=\beta,
\end{equation}
where we used (\ref{eq:find-b-given-a}). This implies that the
left and the right derivatives of $W$ are equal at $F(b)$ (smooth
fit), since the left derivative at $F(b)$ is also equal to
$\beta$.
\end{remark}
\noindent \underline {Third stage: } Now, we vary $a \in
\mathbb{R}_{+}$ and choose $a^{*}$ that maximizes $a \rightarrow
\beta(a)$. We also find the corresponding $b^{*}=b(a^*)$. Now, the
value function is given by (\ref{eq:perform-a-b}) when $a$ and $b$
are replaced by $a^*$ and $b^*$ respectively.

The next proposition justifies the second stage of our algorithm.
\begin{proposition}\label{prop:main-proposition} Assume that $r(a;a)<0$ and $sup_{x} r(x;a)>0$.
Furthermore, if the functions $R(\cdot;a)$ and $H(\cdot)$ defined in
(\ref{eq:transform-r-h}) are increasing and concave on
$(y,\infty)$ for some $y \geq F(a)$, and the function $h(\cdot)$
defined in (\ref{defn:r-h}) satisfies $h(\cdot) \in (0,1)$, then
for any given $a \geq 0$, (\ref{eq:find-b-given-a}) has a unique
solution.
\end{proposition}

The proof essentially follows from Remark~\ref{remark:smooth-fit}.
But we will have to introduce a series of lemmas before we justify
our claim.

First, let us also introduce a family of value functions
parameterized by $\gamma \in \mathbb{R}$ as
\begin{equation}\label{eq:param}
V_a^\gamma(x)\triangleq\sup_{\tau\in \S}\ME\left[e^{-\alpha \tau}
r^{\gamma}(X_{\tau};a)\right] \quad \text{where} \quad r^{\gamma}(x;a)
\triangleq e^{-\alpha \Delta}(r(x;a)+\gamma \cdot h(x)),
\end{equation}
in which $\S$ is the set of stopping times of the natural filtration of $X$. Here, $X$ is a diffusion on $[0,\infty)$, which is absorbed at the left boundary. (In the case of geometric Brownian motion this left boundary is taken to be $d>0$.)
Then we have the following result.
\begin{lemma}\label{lem:majorant-forex}
Let us define
\begin{equation}
R^{\gamma}(\cdot;
a)\triangleq\frac{r^{\gamma}(F^{-1}(\cdot),a)}{\varphi(F^{-1}(\cdot))},
\end{equation}
then the function
\begin{equation}\label{eq:Wa} W^\gamma_a(\cdot) \triangleq
\frac{V^\gamma_a(F^{-1}(\cdot))}{\varphi(F^{-1}(\cdot))},
\end{equation}
is the smallest non-negative concave majorant of $R^{\gamma}$ that
passes through  $(F(0),0)$. Moreover under the assumptions of
Proposition~\ref{prop:main-proposition} this majorant is linear in
the \emph{continuation region} (the region in which $W^{\gamma}$
is strictly greater than $R^{\gamma}$).
\end{lemma}
\begin{proof}
The first part of the proof follows Proposition 5.3 of Dayanik
and Karatzas \cite{DK2003}. The second part of the proof follows
from the first and the fact that $R^{\gamma}(\cdot;a)$ is
increasing and concave on $(y,\infty)$.
\end{proof}

The following technical lemma will be used in showing the
existence of $\gamma$ such that $V_a^{\gamma}(a)=\gamma$ for any
$a\geq 0$.

\begin{lemma}\label{lem:Lipschitz}
\begin{equation}\label{eq:lipschitz-statement}
V_a^{\gamma_1}(x)-V^{\gamma_2}_a(x) \leq \gamma_1-\gamma_2, \quad
\text{for any} \quad \gamma_1 \geq \gamma_2 \quad \text{and} \quad
x \geq 0.
\end{equation}
\end{lemma}
\begin{proof}
It is clear from (\ref{defn:r-h}) that $\gamma \rightarrow
V_a^{\gamma}(x)$ is an increasing convex function. Therefore the
right-derivative
\[
D_{\gamma}^{+}V^{\gamma'}_{a}(x) \triangleq \lim_{h\downarrow
0}\frac{V^{\gamma'+h}(x)-V^{\gamma'}(x)}{h}
\]
exists for any $\gamma'>0$ and it satisfies
\begin{equation}\label{eq:lipschitz}
\frac{V^{\gamma_1}_a(x)-V^{\gamma_2}_a(x)}{\gamma_1-\gamma_2} \leq
D_{\gamma}^{+}V^{\gamma_1}_{a}(x),
\end{equation}
for any $\gamma_1 \geq \gamma_2$ (see e.g. \cite{kn:karat}, pages
213-214). Note that since $h(\cdot) \in (0,1)$ we have that
\begin{equation}\label{eq:right-der-bounded}
0<D_{\gamma}^{+}V^{\gamma'}_{a}(x) \leq 1.
\end{equation}
Now, (\ref{eq:lipschitz}) and (\ref{eq:right-der-bounded})
together imply (\ref{eq:lipschitz-statement}).
\end{proof}

\begin{lemma}\label{lem:existence-lemma}
Under the assumptions of Proposition~\ref{prop:main-proposition},
 there exists a unique $\gamma$ such that
$V_a^\gamma(a)=\gamma$ for $a \geq 0$.
\end{lemma}

\begin{proof}
Consider the function $\gamma \rightarrow V^\gamma_a(a)$. Our aim
is to show that there exists a fixed point to this function. Let
us consider $V^0_a(a)$ first. Since $sup_{x}r(x;a)>0$ we have that
$V^0_a(a)>0$. Now let us consider the case when $\gamma>0$. First,
note that $W^\gamma_a(F(a)) \geq R^{\gamma}(F(a), a)$ for all
$\gamma$. Since by Lemma~\ref{lem:Lipschitz}  $V$ has less than
linear growth in $\gamma$ and $R^{\gamma}$ is linear in $\gamma$,
we can find a $\gamma^{'}$ large enough such that
$W^\gamma_a(F(a))=R^{\gamma}(F(a), a)$ for $\gamma \geq
\gamma^{'}$. This implies however
\begin{align*}
V^{\gamma'}_a(a)=\varphi(a)R^{\gamma'}(F(a);a)=e^{-\alpha
\Delta}(r(a;a)+ \gamma' h(a))<\gamma'.
\end{align*}
Since $\gamma \rightarrow V^{\gamma}_a$ is continuous, which
follows from the fact that this function is convex,  $V^0_a>0$ and
$V^{\gamma^{'}}_a(a)<\gamma^{'}$ implies that $\gamma \rightarrow
V^{\gamma}_a$ crosses the line $\gamma \rightarrow \gamma$. Since
$\gamma \rightarrow V^{\gamma}_a$ is increasing convex it crosses this
line only once.
\end{proof}

\noindent \textbf{Proof of Proposition
~\ref{prop:main-proposition}.} The smallest concave majorant
$W^{\gamma}_a$ in (\ref{eq:Wa}) is linear on $(F(0),F(b^\gamma))$
for a unique $b^{\gamma} \in [0,\infty)$ and smoothly fits to
$R^{\gamma}(\cdot;a)$ at $b^{\gamma}$ and coincides with
$R^{\gamma}(\cdot,a)$ on $[b^{\gamma},d)$. Together with
Lemma~\ref{lem:existence-lemma} this implies that there exists a
unique $\gamma^*$ such that
 equations (\ref{eq:continuous2-fit}) and
(\ref{eq:smooth-fit}) are satisfied when $W$ is replaced by
$W^{\gamma^*}_a$ and $b$ is replaced by $b^{\gamma^*}$. If the
solution of equations (\ref{eq:continuous2-fit}) and
(\ref{eq:smooth-fit}) were not unique, on the other hand, then one
would be able to find multiple smooth fit points $b^{\gamma^*}$,
which yields a contradiction.

\subsection{Are the Assumptions of Proposition~\ref{prop:main-proposition} Satisfied?}\label{sec:satisfaction}

The following remark will be helpful in the analysis that follows:
\begin{remark}\label{rem:concavity}
Given a function $k$ let us denote $K(y) \triangleq \frac{k}{\varphi}\circ F^{-1}(y)$,
$y>0$. If $k$ is twice differentiable at $x \geq 0$ and if we
denote $y \triangleq F(x)$, then $K'(y)=m(x)$ and
$K''(y)=\frac{m'(x)}{F'(x)}$ with
\begin{equation}\label{eq:concavity}
m(x)=\frac{1}{F'(x)}\left(\frac{k}{\varphi}\right)'(x), \quad
\text{and} \quad K''(y)[(\mathcal{A}-\alpha)k(x)] \geq 0, \quad
y=F(x),
\end{equation}
with strict inequality if $H''(y) \neq 0$. The inequality in
(\ref{eq:concavity}) is useful in identifying the concavity of
$K$.
\end{remark}

\subsubsection{Brownian Motion with Drift}

In this case $r(x;a)$ and $h(x)$ defined in (\ref{defn:r-h}) are
given by
\begin{equation}
\begin{split}
r(x;a)&=(x+\mu\Delta-a-\lambda)\Norm+\sigma\sqrt{\Delta}\Density\\
&-e^{-2\mu x/\sigma^2}\left((-x+\mu\Delta - a
-\lambda)\NormNeg+\sigma\sqrt{\Delta}\DensityNeg\right), \\
h(x)&=\left(\Norm - e^{-2\mu x/\sigma^2}\NormNeg\right),
\end{split}
\end{equation}
in which $\phi(x)=(1/\sqrt{2\pi})e^{-x^2/2}$.

First note that $h(x) \in (0,1)$. It is enough to show that
$R(\cdot;a)$ and $H(\cdot)$ are eventually increasing, and are
eventually concave. First, we will show that they are eventually
increasing. The derivative of $R(\cdot;a)$ has the same sign as
\begin{equation}\label{eq:derivative-of-R-BM}
\begin{split}
\varphi(x)^2\left(\frac{r(x;
a)}{\varphi(x)}\right)'(x)&=\varphi(x)\Bigg\{\Norm-\frac{(a+\lambda)}{\sigma\sqrt{\Delta}}\Density\\
&+e^{-2\mu
x/\sigma^2}\left(\NormNeg-\frac{(a+\lambda)}{\sigma\sqrt{\Delta}}\DensityNeg\right)\\
&-D_1{\sigma^2}\left((-x+\mu\Delta-a-\lambda)\NormNeg+\sigma\sqrt{\Delta}\DensityNeg\right)\\
&-D_2e^{-2\mu
x/\sigma^2}\left((x+\mu\Delta-a-\lambda)\Norm+\sigma\sqrt{\Delta}\Density\right)\Bigg\},
\end{split}
\end{equation}
since $F$ is an increasing function. If we take $x>a$ ($a$ is
fixed) large enough, the third line of
(\ref{eq:derivative-of-R-BM}) dominates the the other lines. Since
$D_1>0$, we can conclude that there exists $x' \geq a$ such that
$\left(\frac{r(x; a)}{\varphi(x)}\right)'(x)>0$ on $x\in (x',
\infty)$.

On the other hand, directly taking the derivative, $h(x)$ can be
shown to be an increasing function in $x\in \R_+$, from which it
follows that $H(y)=h(F^{-1}(y))/\varphi(F^{-1}(y))$ is also
increasing.

Next, we will show that $R$ and $H$ are eventually concave.
Consider the equation $(\A-\alpha)r(x; a)=p(x; a)$ so that
\begin{equation*}
p(x; a)=\mu r'(x; a)+\frac{1}{2}\sigma^2 r''(x; a)-\alpha r(x; a).
\end{equation*}
Directly taking the derivatives and letting $x\rightarrow\infty$,
we obtain $r'(x; a)\rightarrow 1$, $r''(x; a)\rightarrow 0$ and
$r(x; a)\rightarrow \infty$. Therefore $\lim_{x\rightarrow
\infty}p(x; a)=-\infty$.    Similarly, we consider the equation
$q(x)\triangleq (\A-\alpha)h(x)$. By letting $x\rightarrow
\infty$, we have $h(x)\rightarrow 1$, $h'(x)\rightarrow 0$ and
$h''(x)\rightarrow 0$ so that $\lim_{x\rightarrow \infty}q(x)<0$.
Together with Remark~\ref{rem:concavity}, these facts imply that
$R(\cdot,a)$ and $H(\cdot)$ are concave on $y\in (y'', +\infty)$
for some $y''F(a)$.

\subsubsection{Ornstein-Uhlenbeck Process}
We will only consider the case when the performance function is as
in (\ref{u-cases}). The analysis for the case when declaring
banktruptcy is penalized can be performed similarly, since first
and the second derivatives of the integral term in
(\ref{eq:tilde-B}) with respect to the $x$ variable  goes to zero
as $x \rightarrow \infty$.

In this case $r(x;a)$ and $h(x)$ defined in (\ref{defn:r-h}) are
given by
\begin{equation}
\begin{split}
r(x;a)&=xe^{-\rho\Delta}-\left(2N\left(\frac{xe^{-\rho\Delta}}{\sqrt{Q(\Delta)}}\right)-1\right)(a+\lambda),
\quad
h(x)=2N\left(\frac{xe^{-\rho\Delta}}{\sqrt{Q(\Delta)}}\right)-1.
\end{split}
\end{equation}
First, observe that $r(x; a)>0$ on $(x_1, \infty)$ with some
$x_1>a$.  By taking the derivative of $r(x; a)$, we have
\begin{align*}
r'(x;
a)=e^{-\rho\Delta}\left(1-\frac{2(a+\lambda)}{\sqrt{Q(\Delta)}}\phi\left(\frac{xe^{-\rho\Delta}}{\sqrt{Q(\Delta)}}\right)\right).
\end{align*}
From this expression we see that $r'(x; a)>0$ on $x\in (x_2,
\infty)$ with some $x_2>a$. Let us denote $x'\triangleq \max(x_1,
x_2)$. It follows that $R(y)$ is increasing on $y\in (y', \infty)$
with $y'=F(x')$ because
\begin{equation*}
\left(\frac{r}{\varphi}\right)'=\frac{r'\varphi-r\varphi'}{\varphi^2}
\quad\text{with}\quad \varphi'<0.
\end{equation*}
Observe also that $h(x) \in (0,1)$ and $h'(x)>0$ on $x\in \R_+$.

Next, we will analyze the concavity properties of $R(\cdot;a)$ and
$H$. Consider the equation $(\A-\alpha)r(x; a)=p(x; a)$ so that
\begin{equation*}
p(x; a)=-\rho xr'(x; a)+\frac{1}{2} r''(x; a)-\alpha r(x; a).
\end{equation*}
We have $r(x; a)\rightarrow +\infty$, $xr'(x; a)\rightarrow
+\infty$ and $r''(x; a)\rightarrow 0$ as $x\rightarrow \infty$.
Thus, we have $\lim_{x\rightarrow \infty}p(x; a)=-\infty$.
Similarly, we consider the equation $q(x)\triangleq
(\A-\alpha)h(x)$. By letting $x\rightarrow \infty$, we have
$h(x)\rightarrow 1$, $xh'(x)\rightarrow 0$ and $h''(x)\rightarrow
0$ so that $\lim_{x\rightarrow \infty}q(x)<0$. Together with
Remark~\ref{rem:concavity}, this analysis shows that there exists
$y'' \geq F(a)$ such that $R(\cdot;a)$ and $H(\cdot)$ are concave
on $(y'',\infty)$.

\subsubsection{Square Root Process}
In this case the functions $r$ and $h$ are given by
\begin{equation}\label{eq:r-h-sqrt}
\begin{split}
r(x;a)&=(xe^{-2\rho\Delta}+Q(\Delta)-(a+\lambda))\left(2\NormSqrt-1\right)+2\sqrt{Q(\Delta)x}e^{-\rho\Delta}\DensitySqrt,
\\ h(x)&=2\NormSqrt-1.
\end{split}
\end{equation}
Observe that $r(x; a)>0$ on $(x_1, \infty)$ with some $x_1 \geq a$
since the only negative term in the first equation in
(\ref{eq:r-h-sqrt}) is bounded from below by $-(a+\lambda)$.
Taking the derivative of $r(x; a)$ we obtain
\begin{equation}\label{eq:der-r-sqrt}
\begin{split}
r'(x;
a)&=e^{-2\rho\Delta}\left(2\NormSqrt-1\right)+\frac{e^{-\rho\Delta}}{\sqrt{Q(\Delta)}}\DensitySqrt\left(\sqrt{x}e^{-2\rho\Delta}+\frac{a+\lambda}{\sqrt{x}}\right)\\
&-\frac{e^{-\rho\Delta}}{\sqrt{x}}\DensitySqrt\left(\sqrt{Q(\Delta)}+\frac{xe^{-\rho\Delta}+Q(\Delta)}{\sqrt{Q(\Delta)}}\right).
\end{split}
\end{equation}
The second term on the first line of (\ref{eq:der-r-sqrt}) is
positive and it dominates as $x \rightarrow \infty$, therefore
$r'(x; a)>0$ on $x\in (x_2, \infty)$ with for some $x_2 \geq a$.
 Take $x'\triangleq \max(x_1, x_2)$. It follows that
$R(y)$ is increasing on $y\in (y', \infty)$, in which $y'=F(x')$.
On the other hand, $h(x) \in (0,1)$ and
$h'(x)=-\frac{e^{-\rho\Delta}}{\sqrt{x}}\DensitySqrt<0$. However,
$h'$ goes to zero as $x \rightarrow \infty$, which implies that
$(h/\varphi)'=\frac{h'\varphi-h\varphi'}{\varphi^2}>0$ on $(x'',
\infty)$ for some sufficiently large $x''$.

Next, we analyze the concavity properties of $R(\cdot;a)$ and 
$H(\cdot)$. Let us define $p(x; a) \triangleq (\A-\alpha)r(x; a)$.
As a result
\begin{equation*}
p(x; a)=(1-2\rho x)r'(x; a)+2x r''(x; a)-\alpha r(x; a).
\end{equation*}
We have $r(x; a)\rightarrow +\infty$, $xr'(x; a)\rightarrow
+\infty$ and $xr''(x; a)\rightarrow 0$ as $x\rightarrow \infty$.
Thus, we have $\lim_{x\rightarrow \infty}p(x; a)=-\infty$.
Similarly, we consider the equation $q(x)\triangleq
(\A-\alpha)h(x)$. By letting $x\rightarrow \infty$, we have
$h(x)\rightarrow 1$, $xh(x)\rightarrow 0$ and $xh''(x)\rightarrow
0$ so that $\lim_{x\rightarrow \infty}q(x)<0$. Using
Remark~\ref{rem:concavity}, we observe that $R(\cdot;a)$ and
$H(\cdot)$ are eventually concave.

\subsubsection{Geometric Brownian Motion}

When the aggregate income process $X^0$ is modeled by a geometric
Brownian motion a sufficient condition for the hypothesis of the
Proposition~\ref{prop:main-proposition} to hold is $\mu \leq
\alpha$. In this case the functions $r$ and $h$ are given by
\begin{equation}\label{eq:r-Gbm}
\begin{split}
r(x;a)&=e^{\mu\Delta}x\left(N(d_1)-\left(\frac{d}{x}\right)^{1+2\mu/\sigma^2}N(-d_2)\right)
-(a+\lambda)\left(N\left(\frac{-\tilde{d}+\gamma\Delta}{\sigma\sqrt{\Delta}}\right)-e^{2\gamma\tilde{d}/\sigma^2}
N\left(\frac{\tilde{d}+\gamma\Delta}{\sigma\sqrt{\Delta}}\right)\right),
\\ h(x)&=\left(N\left(\frac{-\tilde{d}+\gamma\Delta}{\sigma\sqrt{\Delta}}\right)-e^{2\gamma\tilde{d}/\sigma^2}
N\left(\frac{\tilde{d}+\gamma\Delta}{\sigma\sqrt{\Delta}}\right)\right).
\end{split}
\end{equation}
Observe that $\tilde{d}<0$  since $x>d$. Moreover,
$N(d_1)\rightarrow 1$, $N(-d_2)\rightarrow 0$ and
$N(\tilde{d})\rightarrow 0$ as $x\rightarrow +\infty$. Also, $r(x;
a)>0$ on $(x_1, \infty)$ with some $x_1>a$ since the negative term
in the first equation in (\ref{eq:r-Gbm}) is bounded.  On the
other hand $h(x)>0$ for $x\in\R_+$ and $h'(x)\rightarrow 0$ as
$x\rightarrow \infty$. The derivative of $r$ is
\begin{align*}
r'(x;
a)&=e^{\mu\Delta}\left(N(d_1)-\left(\frac{d}{x}\right)^{1+2\mu/\sigma^2}N(-d_2)\right)-(a+\lambda)h'(x)
\\
&+e^{\mu\Delta}\left(\frac{\phi(d_1)}{\sigma\sqrt{\Delta}}+(1+2\mu/\sigma^2)\left(\frac{d}{x}\right)^{1+2\mu/\sigma^2}
+\left(\frac{\phi(-d_2)}{\sigma\sqrt{\Delta}}\right)\left(\frac{d}{x}\right)^{1+2\mu/\sigma^2}\right),
\end{align*}
which is positive on $x\in (x_2, \infty)$ for some $x_2 \geq a$.
Take $x'\triangleq \max(x_1, x_2)$. It follows that $R(y)$ is
increasing on $y\in (y', \infty)$ with $y'=F(x')$.  Similarly,
since $h(x) \in (0,1)$ and $h'(x)$ goes to zero as $x \rightarrow
\infty$, so that $H'(y)>0$ on $(y'', \infty)$ for sufficiently
large $y''$.

Next, we analyze the concavity of $R(\cdot;a)$ and $H(\cdot)$. Let
us denote $p(x; a) \triangleq (\A-\alpha)r(x; a)$. The function
$p(\cdot;a)$ is given by
\begin{align*}
p(x; a)&=\mu xr'(x; a)+\frac{1}{2} \sigma^2 r''(x; a)-\alpha r(x;
a)\\
&=(\mu-\alpha)xe^{\mu\Delta}\left(N(d_1)-\left(\frac{d}{x}\right)^{1+2\mu/\sigma^2}N(-d_2)\right)
-\alpha(a+\lambda)h(x)+T(x; a),
\end{align*}
where $T(x; a)$ is the terms that involve $\phi(\cdot)$ or
$\left(\frac{d}{x}\right)^{1+2\mu/\sigma^2}$ and
$\lim_{x\rightarrow +\infty}T(x; a)=0$.  Observe that
$\lim_{x\rightarrow +\infty}p(x; a)=-\infty$ when $\mu \leq
\alpha$. Similarly, we consider the equation $q(x)\triangleq
(\A-\alpha)h(x)$. $h(x)\rightarrow 1$, $h'(x)\rightarrow 0$ and
$h''(x)\rightarrow 0$ implies that $\lim_{x\rightarrow
\infty}q(x)<0$. Using Remark~\ref{rem:concavity}, we can conclude
that $R(\cdot;a)$ and $H(\cdot)$ are eventually concave.

\section{Numerical Examples}

See Figures 1-4 for numerical illustrations. In our examples we
quantify the effect of delay in dividend payments. In each case we
find the optimal dividend payment barrier, $b^*$ , the optimal
amount of dividend payment, $b^{*}-a^*$, and the value function
$v$. Then we compare them to $b^0$, $b^0-a^0$ and $v^0$, the
analogues of the previous quantities when there is no delay. As
expected the value function is smaller, $v<v_0$ when there is
delay in dividend payments. Since in Figures 2 (b), 2 (e) and 4
(b), the value functions $v$ and $v_0$ are not distinguishable, in
Figures (2) and (4) we plot the difference of $v_0-v$.

When the aggregate income process, $X^0$ is modeled by a Brownian
motion with drift, a square root process then we observe that
$a^*<a_0$, $b^*<b_0$, $b^*-a^*<b_0-a_0$ and $\beta^*<\beta_0$. The same conclusion
holds if  $X^0$ is an Ornstein-Uhlenbeck process and the declaring
bankruptcy is penalized. On the other hand, when $X^0$ is modeled
by an Ornstein-Uhlenbeck process (the case in which declaring ruin
is not penalized) or a geometric Brownian motion we obtain that
$a^*=a_0$, $b^*>b_0$, $b^*-a^*>b_0-a_0$ and $\beta^*<\beta_0$. Note that in both of
these cases declaring bankruptcy is optimal as soon as the
aggregate income level hits $b^*$, regardless of the magnitude of
delay.

Observe that in the numerical examples considered, the function
$\beta(a)$, which is obtained from (\ref{eq:slope}) after we plug
in for $b$ that we obtain from (\ref{eq:find-b-given-a}) (say
$b(a)$), is concave. It is either strictly decreasing or has a
unique local maximum. We leave the proof of these features of the
function $\beta(a)$ as an open problem.

\begin{remark} \normalfont In our framework, it is easy to deal with solvency
constraints. The optimal $a^*$ may not be acceptable,  and
prohibited by regulatory constraints. This was studied by Paulsen
\cite{paulsen} in singular control setting (with no delays). Let
as consider the case with $\Delta=0$ and assume that the firm is
not allowed to reduce its aggregate cash flow to below
$\tilde{a}$. If we show the above properties hold for $\beta(a)$
$a$ it is easy to argue if $a^*>\tilde{a}$, then every time it
pays out dividends the firm would reduce its reservoir to $a^*$
(the constraint is not binding), else if $a^{*}< \tilde{a}$, then
the firm every time it pays out dividends the firm would reduce
its reservoir to $\tilde{a}$.
\end{remark}

\begin{figure}\label{fig:BM}[h]
\begin{center}
\begin{minipage}{0.45\textwidth}
\centering \includegraphics[scale=0.85]{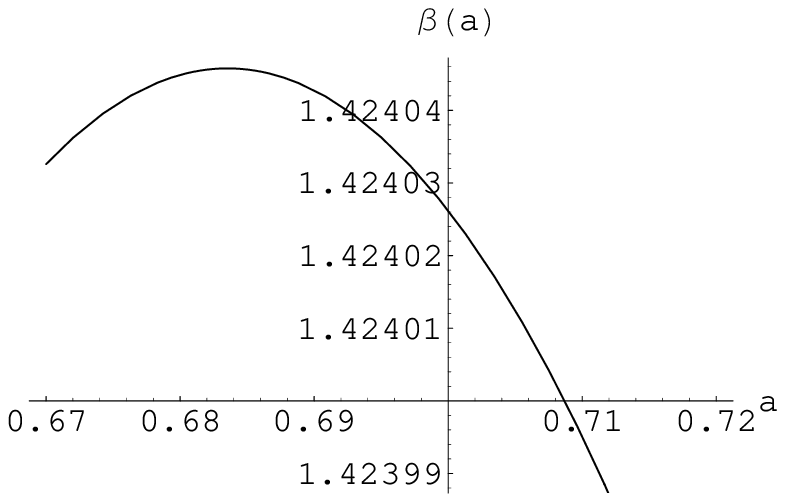} \\ (a)
\end{minipage}
\begin{minipage}{0.45\textwidth}
\centering \includegraphics[scale=0.85]{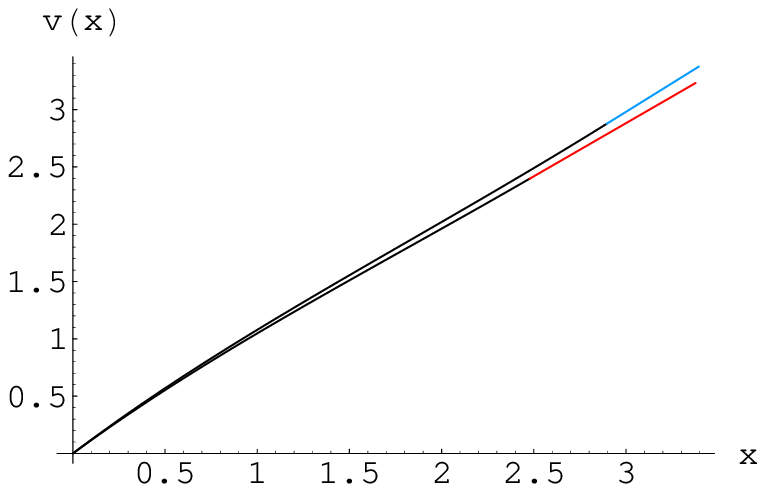} \\ (b)
\end{minipage}
\caption{\small A numerical example of a Brownian motion with
drift with parameters $(\mu, \alpha, \sigma, \lambda,
\Delta)=(0.3, 0.15, \sqrt{2}, 0.1, 0.25)$: (a) The graph of
$\beta(a)$ that attains the global maximum at $a^*=0.755$ with
$\beta^*=1.443$. (b) The value function $v(x)$ (below) with
$b^*=2.719$. It is compared with the case of $\Delta=0$ (above)
with $(a_0, b_0, \beta_0)$=$ (0.850, 2.895, 1.466)$. }
\end{center}
\end{figure}

\begin{figure}\label{fig:OU}[h]
\begin{center}
\begin{minipage}{0.45\textwidth}
\centering \includegraphics[scale=0.85]{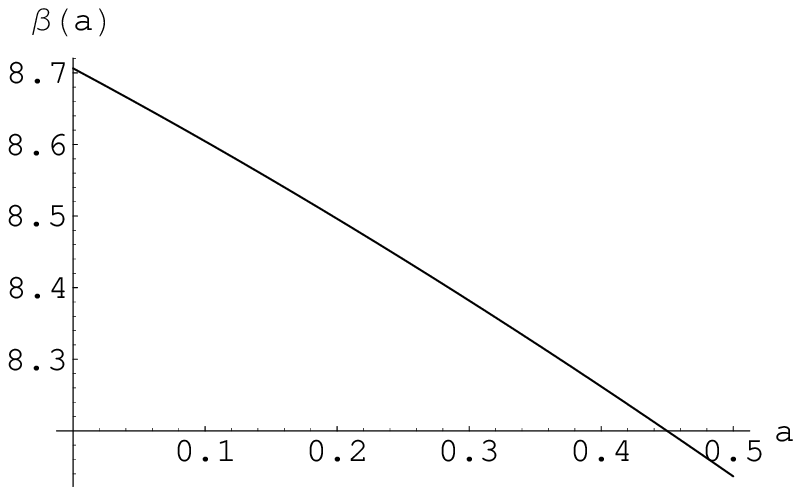} \\ (a)
\end{minipage}
\begin{minipage}{0.45\textwidth}
\centering \includegraphics[scale=0.85]{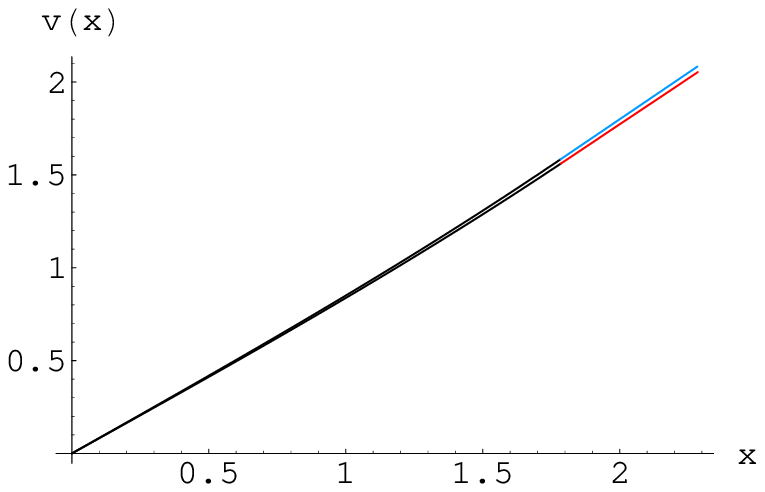} \\ (b)
\end{minipage}
\begin{minipage}{0.45\textwidth}
\centering \includegraphics[scale=0.85]{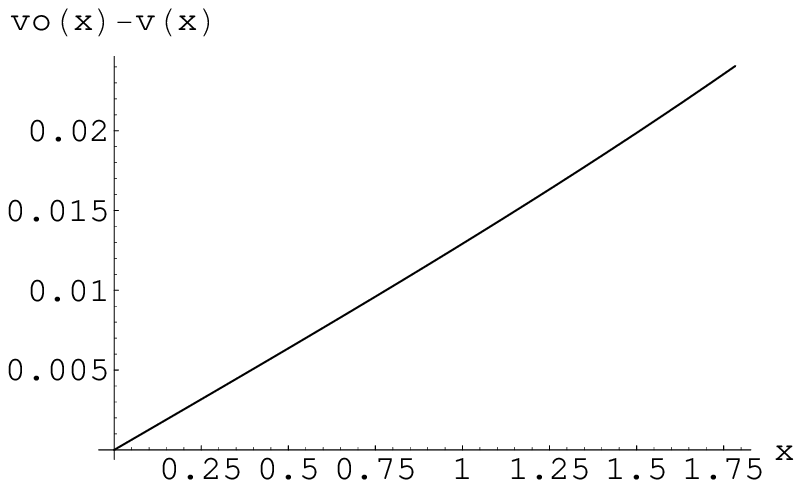} \\ (c)
\end{minipage}
\begin{minipage}{0.45\textwidth}
\centering \includegraphics[scale=0.85]{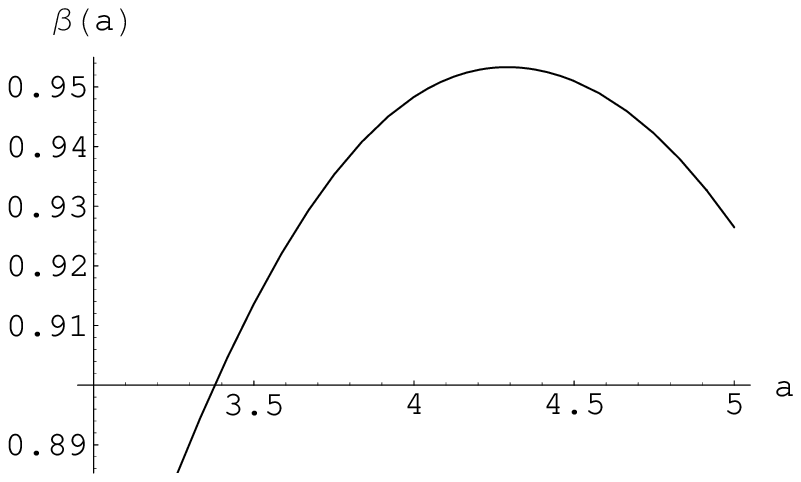} \\ (d)
\end{minipage}
\begin{minipage}{0.45\textwidth}
\centering \includegraphics[scale=0.85]{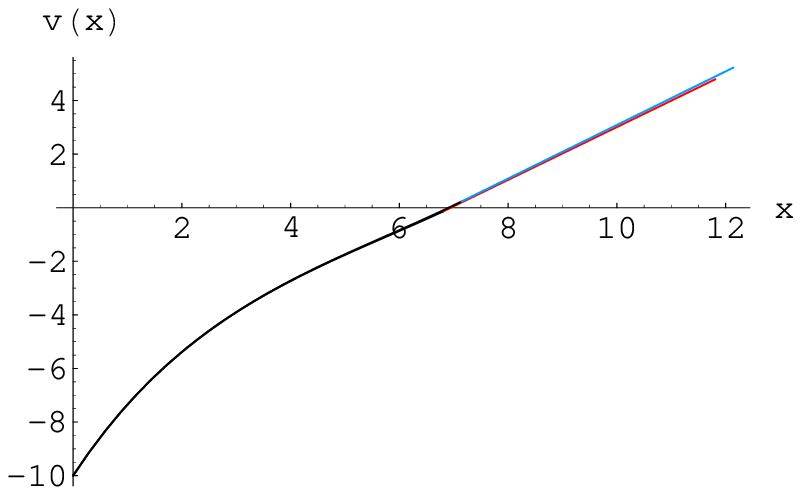} \\ (e)
\end{minipage}
\begin{minipage}{0.45\textwidth}
\centering \includegraphics[scale=0.85]{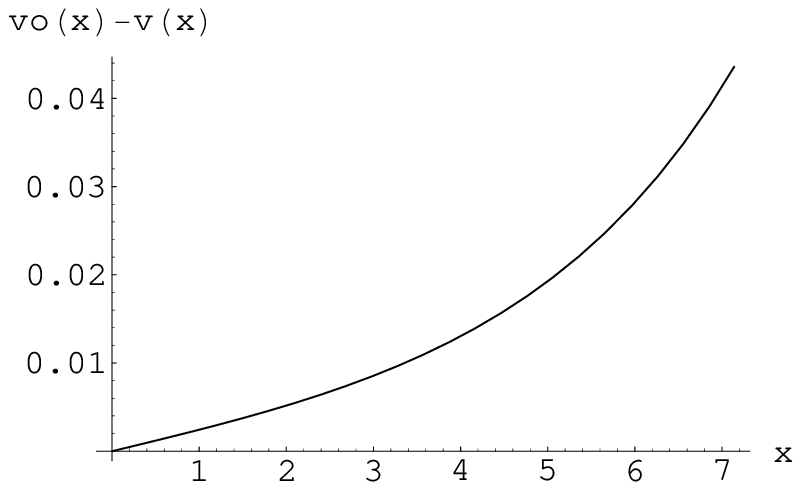} \\ (f)
\end{minipage}
\caption{\small A numerical example of an OU process with
parameters $(\rho, \alpha,  \lambda, \Delta)=(0.01, 0.05,
0.2,0.25)$: (a) The graph of $\beta(a)$ that attains the global
maximum at $a^*=0$ with $\beta^*=8.706$. (b) The value function
$v(x)$ (below) with $b^*=1.785$. It is compared with the case
$v_0(x)$ of $\Delta=0$ (above) with $(a_0, b_0, \beta_0)$=$ (0,
1.783, 8.841)$.  (c) Plot of the  difference  $v_0(x)-v(x)$. (d)
In the case of penalty at ruin, $P=10$, the graph of $\beta(a)$
that attains the global maximum at $a^*=4.290$ with
$\beta^*=0.953$.  (b) The value function $v(x)$ (below) with
$b^*=6.811$.  It is compared with the case $v_0(x)$ of $\Delta=0$
(above) with $(a_0, b_0, \beta_0)=(4.349, 7.141, 0.979)$.  (e)
Plot of the difference $v_0(x)-v(x)$.}
\end{center}
\end{figure}

\begin{figure}[h]
\begin{center}\label{fig:SR}
\begin{minipage}{0.45\textwidth}
\centering \includegraphics[scale=0.85]{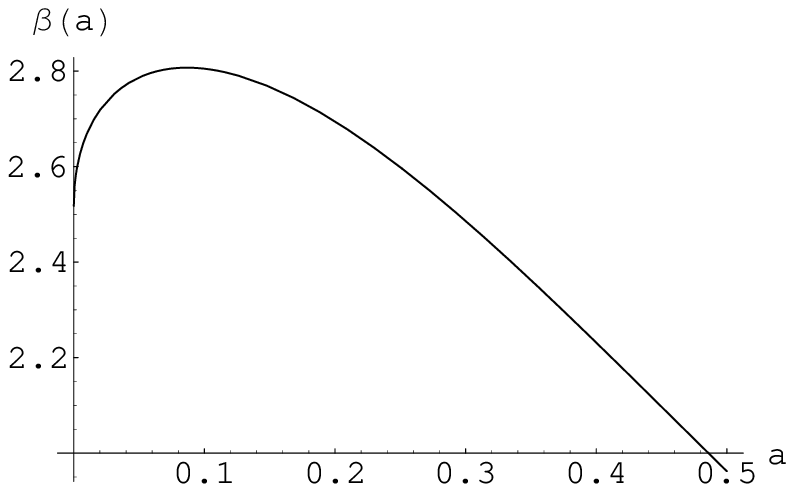} \\ (a)
\end{minipage}
\begin{minipage}{0.45\textwidth}
\centering \includegraphics[scale=0.85]{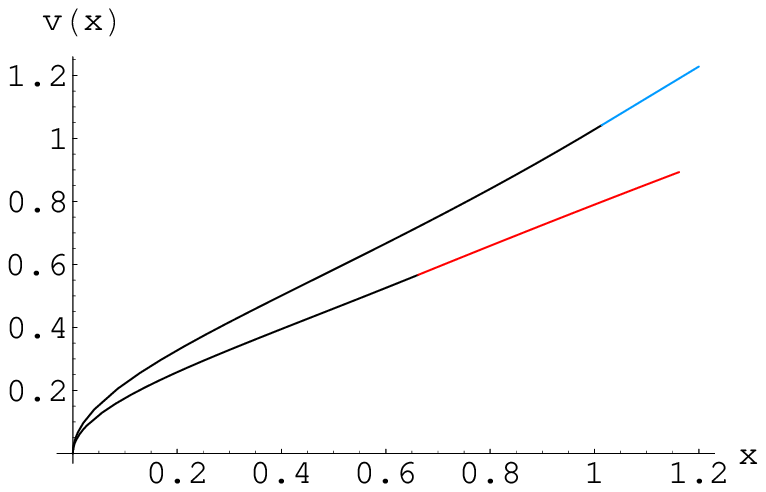} \\ (b)
\end{minipage}
\caption{\small A numerical example of a square root  process with
parameters $(\rho, \alpha,  \lambda, \Delta)=(1, 0.1, 0.1, 0.25)$:
(a) The graph of $\beta(a)$ that attains the global maximum at
$a^*=0.09$ with $\beta^*=2.807$. ( b) The value function $v(x)$
(below) with $b^*=0.662$. It is compared with the case of
$\Delta=0$ (above) with $(a_0, b_0, \beta_0)$=$ (0.165, 1.014,
3.561 )$. }
\end{center}
\end{figure}

\begin{figure}\label{fig:GB}[h]
\begin{center}
\begin{minipage}{0.45\textwidth}
\centering \includegraphics[scale=0.85]{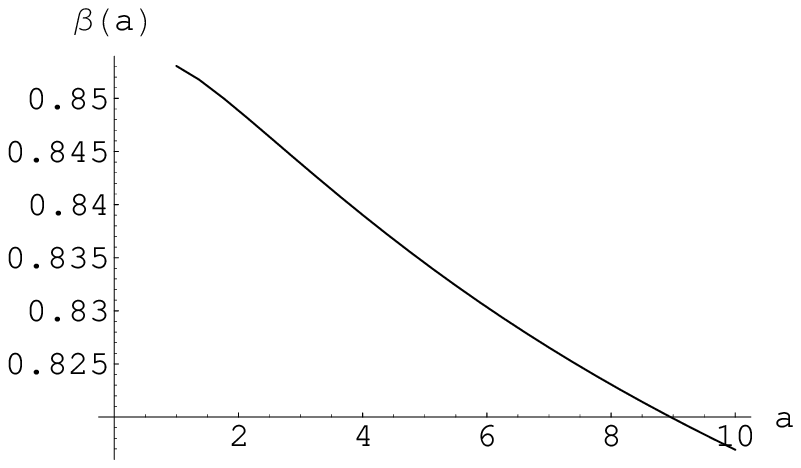} \\ (a)
\end{minipage}
\begin{minipage}{0.45\textwidth}
\centering \includegraphics[scale=0.85]{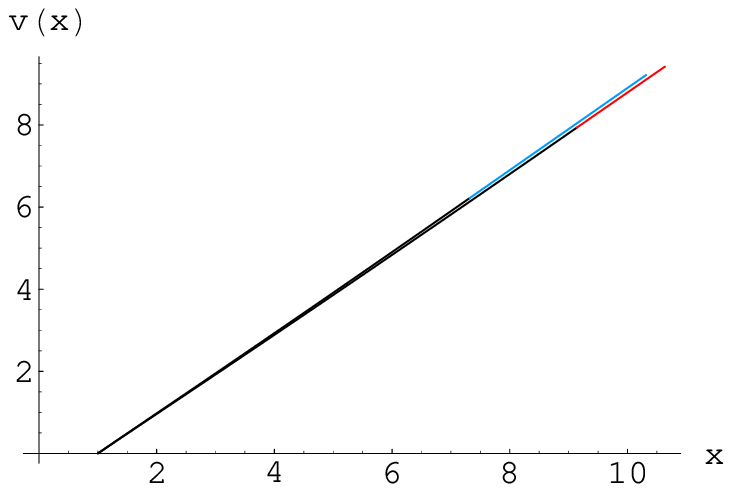} \\ (b)
\end{minipage}
\begin{minipage}{0.45\textwidth}
\centering \includegraphics[scale=0.85]{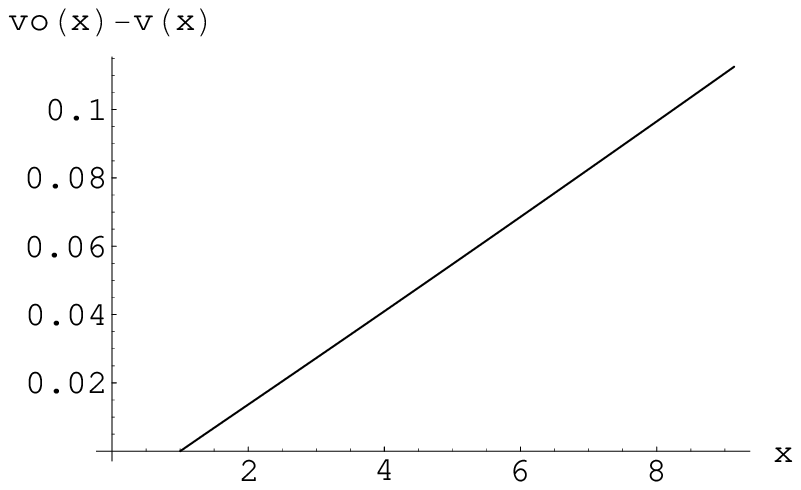} \\ (c)
\end{minipage}
\caption{\small A numerical example of an geometric Brownian
motion with parameters $(\mu, \sigma, \alpha,  \lambda,
\Delta)=(0.05, \sqrt{2}, 0.1, 0.1, 0.25)$ and the ruin level
$d=1$: (a) The graph of $\beta(a)$ that attains the global maximum
at $a^*=1$ with $\beta^*=0.853$. (b) The value function $v(x)$
(below) with $b^*=9.138$. It is compared with the case of
$\Delta=0$ (above) with $(a_0, b_0, \beta_0)$=$ (1, 7.318,
0.865)$.  (c) Plot of the difference $v_0(x)-v(x)$.}
\end{center}
\end{figure}

\section{Conclusion}
We study optimal dividend payout problems with delay using various types of diffusions.  Our method facilitates greatly the solution procedure due to the new characterization of the value function.  The existence of the finite value function and the uniqueness of optimal threshold strategy reduce to verifications of the assumption of Proposition \ref{prop:main-proposition}.  Our models here are more realistic since the delays with respect to dividend payments are explicitly handled.

{\small
\bibliographystyle{plain}

}

\end{document}